\documentclass[12pt]{amsart}

\usepackage{amssymb}

\usepackage{enumerate}

\usepackage{graphicx}

\makeatletter
\@namedef{subjclassname@2010}{%
  \textup{2010} Mathematics Subject Classification}
\makeatother

\newtheorem{theorem}{ Main Theorem}[section]

\theoremstyle{definition}

\numberwithin{equation}{section}

\begin{document}


\baselineskip=17pt


\title[On the Diophantine equations  $ \sum_{i=1}^n a_ix_{i} ^6+\sum_{i=1}^m b_iy_{i} ^3= \sum_{i=1}^na_iX_{i}^6\pm\sum_{i=1}^m b_iY_{i} ^3 $ ]{On the Diophantine equations $ \sum_{i=1}^n a_ix_{i} ^6+\sum_{i=1}^m b_iy_{i} ^3= \sum_{i=1}^na_i X_{i}^6\pm\sum_{i=1}^m b_iY_{i} ^3 $  }

\author[F. Izadi]{Farzali Izadi}
\address{Farzali Izadi \\
Department of Mathematics \\ Faculty of Science \\ Urmia University \\ Urmia 165-57153, Iran}
\email{f.izadi@urmia.ac.ir}

\author[M. Baghalaghdam]{Mehdi Baghalaghdam}
\address{Mehdi Baghalaghdam \\
Department of Mathematics\\ Faculty of Science \\ Azarbaijan Shahid Madani University\\Tabriz 53751-71379, Iran}
\email{mehdi.baghalaghdam@azaruniv.edu}

\date{}

\begin{abstract}
In this paper, the elliptic curves theory is used for solving the Diophantine equations
$\sum_{i=1}^n a_ix_{i} ^6+\sum_{i=1}^m b_iy_{i} ^3= \sum_{i=1}^na_iX_{i}^6\pm\sum_{i=1}^m b_iY_{i} ^3$,
where $n$, $m$ $\geq 1$ and $a_i$, $b_i$, are fixed arbitrary nonzero integers. By our method, we may find infinitely many nontrivial positive solutions and also obtain infinitely many nontrivial parametric solutions for the Diophantine equations for every arbitrary integers $n$, $m$, $a_i$ and $b_i$.
\end{abstract}



\maketitle

\section{Introduction}
Number theory is a vast and fascinating field of mathematics, sometimes called ''higher arithmetic'', consisting
of the study of the properties of whole numbers. In mathematics, a Diophantine equations  is a polynomial equation, usually in two or more unknowns, such that only the integer solutions are studied. The word Diophantine refers to the Hellenistic mathematician of the 3rd century, Diophantus of Alexanderia, who made a study of such equations. While individual equation present  a kind of puzzle and have been considered throughout  the history, the formulation of general theories of Diophantine equations (beyond the theory of quadratic forms) was an achievement of the twentieth century. As a history work, Euler conjectured in $1969$ that the Diophantine equation $A^4+B^4+C^4=D^4$, or more generally $A_1^N+A_2^N+ \cdots +A_{N-1}^N=A_N^N$, ($N\geq4$), has no solution in positive integers ( see \cite{1}). Nearly two centuries later, a computer search (see \cite{3})  found the first counterexample to the general conjecture (for $N=5$):
$27^5+84^5+110^5+133^5=144^5$.

 In $1986$, Noam Elkies found a method to construct an infinite series counterexamples for the $K=4$ case (see \cite{2}). His smallest counterexample was:

 $2682440^4+15365639^4+18796760^4=20615673^4$.
\\

In this paper, we are interested in the study of Diophantine equations:

 $ \sum_{i=1}^n a_ix_{i} ^6+\sum_{i=1}^m b_iy_{i} ^3= \sum_{i=1}^na_iX_{i}^6\pm\sum_{i=1}^m b_iY_{i} ^3 $,

where $n$, $m$ $\geq 1$ and $ a_i$, $b_i$, are  fixed arbitrary nonzero integers.

Our main results are the following theorems:

\begin{theorem} Consider the Diophantine equation :
\\

$ \sum_{i=1}^n a_ix_{i} ^6+\sum_{i=1}^m b_iy_{i} ^3= \sum_{i=1}^na_iX_{i}^6+\sum_{i=1}^m b_iY_{i} ^3 $,
\\

where $n$, $m$ $\geq 1$  and $a_i$, $b_i$, are fixed arbitrary nonzero integers. Let $Y^2=X^3+FX^2+GX+H$, be an elliptic curve in which the coefficients $F$, $G$, and $H$ are all functions
of $a_i$, $b_i$, and the other rational parameters  $B_i$, $C_i$, $D_i$, and $h$, yet to be found later. If the elliptic curve has positive rank, depending on the values of $B_i$, $C_i$, $D_i$, and $h$,  the Diophantine equation has infinitely many  integer solutions.
\\

Proof. Let: $x_i=-D_iu+h$, $y_i=B_iu+C_iv$, $X_i=D_iu+h$, and

$Y_i=-B_iu+C_iv$, where all variables are rational numbers. By substituting these variables in the above Diophantine equation, and after some simplifications, we get:
 \\

\begin{equation}\label{8}
v^2=\frac{(2h\sum_{i=1}^n a_iD_i^5)}{(\sum_{i=1}^m b_iB_iC_i^2)} u^4+\frac{(20h^3\sum_{i=1}^n a_iD_i^3- \sum_{i=1}^m b_iB_i^3)}{(3\sum_{i=1}^m b_iB_iC_i^2)}u^2+\frac{(2h^5\sum_{i=1}^n a_iD_i)
}{(\sum_{i=1}^m b_iB_iC_i^2)}.
\end{equation}
\\

To complete the proof we use the following theorem for transforming this quartic to an elliptic curve of the form
$y^2+a_1xy+a_3y=x^3+a_2x^2+a_4x+a_6$, where $a_i\in \mathbb{Q}$.  (see \cite{5})
\\

THEOREM. Let K be a field of characteristic not equal to $2$. Consider the equation

$v^2=au^4+bu^3+cu^2+du+q^2$, with $a$, $b$, $c$, $d$ $\in K$.
\\

Let $x=\frac{2q(v+q)+du}{u^2}$, $y=\frac{4q^2(v+q)+2q(du+cu^2)-(\frac{d^2u^2}{2q})}{u^3}$.
\\

Define
$a_1=\frac{d}{q}$, $a_2=c-(\frac{d^2}{4q^2})$, $a_3=2qb$, $a_4=-4q^2a$, $a_6=a_2a_4$.
\\

Then $y^2+a_1xy+a_3y=x^3+a_2x^2+a_4x+a_6$.
\\

The inverse transformation is
$u=\frac{2q(x+c)-(\frac{d^2}{2q})}{y}$, $v=-q+\frac{u(ux-d)}{2q}$.
\\

The point $(u, v)=(0,q)$ corresponds to the point $(x, y)=\infty$ and

$(u, v)=(0,−q)$ corresponds to $(x, y)=(−a_2,a_1a_2-a_3)$.
\\

We see that if
\begin{equation*}
L:=\frac{(2h^5\sum_{i=1}^n a_iD_i)
}{(\sum_{i=1}^m b_iB_iC_i^2)}
\end{equation*}
\\

 to be square, say $q^2$, (This is possible by fixing the other parameters and  choosing appropriate values for $h$.),
we may use the above theorem and transform the quartic elliptic curve to an elliptic curve of  the form

 $y^2+a_1xy+a_3y=x^3+a_2x^2+a_4x+a_6$, where $a_i\in \mathbb{Q}$.
\\

If the above elliptic curve has positive rank (For every arbitrary nonzero integers $a_i$ and $b_i$, this is done  by choosing appropriate values for $B_i$, $C_i$, $D_i$, $h$.), by calculating $x_i$, $y_i$, $X_i$, $Y_i$, from relations $x_i=-D_iu+h$, $y_i=B_iu+C_iv$,
$X_i=D_iu+h$, and $Y_i=-B_iu+C_iv$, after  some simplifications  and canceling  the  denominators of $x_i=-D_iu+h$, $y_i=B_iu+C_iv$, $X_i=D_iu+h$, and $Y_i=-B_iu+C_iv$, we may obtain infinitely many integer solutions for the  Diophantine equation. The proof is complete.
\end{theorem}

Now we are going to solve some couple of examples:
\\

Example 1. We wish to solve the Diophantine equation:

 $x_1^6+y_1^3=X_1^6+Y_1^3$,
\\

we may assume that in the corresponding the elliptic curve \eqref{8}, $B_1=C_1=D_1=1$,

$h=2$. This means that we let

$x_1=-u+2$, $y_1=u+v$, $X_1=u+2$, $Y_1=-u+v$.
\\

( We may choose  the other appropriate values for $B_1$, $C_1$, $D_1$, and $h$, so that the rank of the corresponding elliptic curve  to be $\geq1$ and $L$ to be square.

Then we get the elliptic curve:

 $v^2=4u^4+53u^2+64$.
\\

With the inverse transformation
 $u=\frac{16(X+53)}{Y}$ and $v=-8+\frac{u^2X}{16}$, this maps to the new elliptic curve

 $Y^2=X^3+53X^2-1024X-54272$.

The rank of this elliptic curve is $1$ and its generator is the point

 $P=(X,Y)=(-48,80)$. Because of this,  the  above elliptic curve has infinitely many rational points and we may obtain infinitely many solutions for the Diophantine equation too. Since $P=(X,Y)=(-48,80)$, we get $(u,v)=(1,-11)$, by calculating $x_1$, $y_1$, $X_1$, $Y_1$, from the above relations and after  some simplifications  and canceling  the  denominators of  $x_1$, $y_1$, $X_1$, $Y_1$, we get the identity:
 \\

$1^6+12^3=3^6+10^3$.
\\

Also we have: $2P=(X',Y')=(68,-660)$,
\\

 $3P=
(X'',Y'')=(\frac{-27168}{841},\frac{490880}{24389})$ and  $4P=(X''',Y''')=(\frac{1139524}{27225},\frac{1180380068}{4492125})$.
\\

By using these three new points, we obtain the other solutions for the Diophantine equation, respectively as 
\\

$74^6+5768^3=14^6+7088^3$,
\\

$1503^6+6175210^3=1919^6+5819322^3$,
\\

$693358^6+2867762572208^3=1427882^6+2478237149768^3$.
\\

By choosing the other points on the elliptic curve such as

$nP$ ($n=5, 6$, $\cdots$ ) we obtain infinitely many solutions for the  Diophantine equation.
\\

Example 2. Let us solve the Diophantine equation:

 $2x_1^6+y_1^3=2X_1^6+Y_1^3$,
\\

we may assume that in the elliptic curve \eqref{8},

$B_1=3$, $D_1=C_1=1$, $h=\frac{4}{3}$. This means that we let

$x_1=-u+\frac{4}{3}$, $y_1=3u+v$, $X_1=u+\frac{4}{3}$, $Y_1=-3u+v$.
\\

( We may choose  the other appropriate values for $B_1$, $C_1$, $D_1$, and $h$, so that rank of the corresponding elliptic curve  to be $\geq1$ and $L$ to be square.

Then we get the elliptic curve:

 $v^2=\frac{16}{9}u^4+\frac{1831}{243}u^2+(\frac{64}{27})^2$.
\\

With the inverse transformation
 $u=\frac{\frac{128}{27}(X+\frac{1831}{243})}{Y}$ and $v=-\frac{64}{27}+\frac{u^2X}{\frac{128}{27}}$, we get the cubic elliptic curve

 $Y^2=X^3+\frac{1831}{243}X^2-\frac{262144}{6561}X-\frac{479985664}{1594323}$.
\\

The rank of this elliptic curve is $2$ and its generators are the points
$P_1=(X,Y)=(\frac{-1792}{243},\frac{3328}{2187})$ and $P_2=(X',Y')=(\frac{-2411}{324},\frac{7007}{5832})$. Because of this,  the  above elliptic curve has infinitely many rational points and we may obtain infinitely many solutions for the Diophantine equation too. By using  this point $P_1=(X,Y)=(\frac{-1792}{243},\frac{3328}{2187})$, we get $(u,v)=(\frac{1}{2},\frac{-149}{54})$,  and by calculating $x_1$, $y_1$, $X_1$, $Y_1$, from the above relations and after some simplifications  and canceling  the  denominators of  $x_1$, $y_1$, $X_1$, $Y_1$, we get the identity:
 \\

$2.(15)^6+1380^3=2.(33)^6+408^3$.
\\

Also we have: $2P_1=(\frac{16708}{729},\frac{-2392940}{19683})$.
\\

By using  these  two new  points $2P_1$, $P_2$, we obtain the other solutions for the Diophantine equation, respectively as
\\

$2.(22773)^6+72216552^3=2.(1317)^6+653700972^3$,
\\

$2.(501)^6+997572^3=2.(885)^6+398820^3$.
\\

By choosing the other points on the elliptic curve such as

$nP_1$ or $nP_2$ ($n=3, 4$, $\cdots$ ) we obtain infinitely many solutions for the  Diophantine equation.
\\

Example 3. We wish to solve the Diophantine equation:

$x_1^6+7y_1^3=X_1^6+7Y_1^3$,
\\

we may assume that in the elliptic curve \eqref{8},

 $B_1=C_1=D_1=1$, $h=\frac{2}{7}$.

 This means that we let $x_1=-u+\frac{2}{7}$, $y_1=u+v$, $X_1=u+\frac{2}{7}$,
$Y_1=-u+v$.
\\

Then we get the elliptic curve:

 $v^2=\frac{4}{49}u^4-\frac{747}{2401}u^2+(\frac{8}{343})^2$.
\\

With the inverse transformation
 $u=\frac{\frac{16}{343}(X-\frac{747}{2401})}{Y}$ and $v=-\frac{8}{343}+\frac{u^2X}{\frac{16}{343}}$,

the cubic elliptic curve is 

 $Y^2=X^3-\frac{747}{2401}X^2-\frac{1024}{5764801}X+\frac{764928}{13841287201}$.
\\

The rank of this elliptic curve is $1$ and its generator is the point

$P=(X,Y)=(\frac{752}{2401},\frac{240}{16807})$. Because of this,  the  above elliptic curve has infinitely many rational points and we may obtain infinitely many solutions for the Diophantine equation too. By using this point $P=(X,Y)=(\frac{752}{2401},\frac{240}{16807})$, we get $(u,v)=(\frac{1}{147},\frac{-3481}{151263})$,  and by calculating $x_1$, $y_1$, $X_1$, $Y_1$, from the above relations and after some simplifications  and canceling  the  denominators of  $x_1$, $y_1$, $X_1$, $Y_1$, we get the identity:
 \\

$287^6+7.(31570)^3=301^6+7.(17164)^3$.
\\

 If we take $B_1=C_1=D_1=1$, $h=\frac{7}{2}$,

 this means that in the above Diphantine equation, we have $x_1=-u+\frac{7}{2}$,

$y_1=u+v$, $X_1=u+\frac{7}{2}$,
$Y_1=-u+v$,
\\

this in turn results to

 $v^2=u^4-\frac{81}{2}u^2+(\frac{49}{4})^2$.
\\

With the inverse transformation
 $u=\frac{\frac{49}{2}(X+\frac{81}{2})}{Y}$ and $v=\frac{-49}{4}+\frac{u^2X}{\frac{49}{2}}$,

we get the new elliptic curve

 $Y^2=X^3+\frac{81}{2}X^2-\frac{2401}{4}X-\frac{194481}{8}$.
\\

The rank of this elliptic curve is $1$ and its generator is the point

$P=(X,Y)=(\frac{-77}{2},42)$. Because of this,  the  above elliptic curve has infinitely many rational points and we may obtain infinitely many solutions for the Diophantine equation too. By using this point $P=(X,Y)=(\frac{-77}{2},42)$, we get $(u,v)=(\frac{7}{6},\frac{-259}{18})$,  and by calculating $x_1$, $y_1$, $X_1$, $Y_1$, from the above relations and after  some simplifications  and canceling  the  denominators of  $x_1$, $y_1$, $X_1$, $Y_1$, we get the identity:
 \\

$7^6+7.(140)^3=14^6+7.(119)^3$.
\\

 Also we have: $2P=(X',Y')=(\frac{2009}{18},\frac{-36260}{27})$.
\\

By using this new  point $2P$, we obtain another solution for the Diophantine equation:
\\

 $251^6+7.(32089)^3=29^6+7.(40969)^3$
 \\

By choosing the other points on the elliptic curve such as

$nP$ ($n=3, 4$, $\cdots$ ) we obtain infinitely many solutions for the  Diophantine equation.
\\
\begin{theorem} Consider the Diophantine equation:
\\

$ \sum_{i=1}^n a_ix_{i} ^6+\sum_{i=1}^m b_iy_{i} ^3= \sum_{i=1}^na_iX_{i}^6-\sum_{i=1}^m b_iY_{i} ^3 $,
\\

where $n$, $m$ $\geq 1$  and $a_i$, $b_i$, are fixed arbitrary nonzero integers. Let $Y^2=X^3+FX^2+GX+H$, be an elliptic curve in which the coefficients $F$, $G$, and $H$ are all functions
of $a_i$, $b_i$, and the other rational parameters  $B_i$, $C_i$, $D_i$, and $h$, yet to be found later. If the elliptic curve has positive rank, depending on the values of $B_i$, $C_i$, $D_i$, and $h$,  the Diophantine equation has infinitely many  integer solutions.
\\

Proof. Let: $x_i=D_iu+h$, $y_i=B_iu-C_iv$, $X_i=-D_iu+h$, and

$Y_i=B_iu+C_iv$, where all variables are rational numbers. By substituting these variables in the above Diophantine equation, and after some simplifications, we get:
 \\

\begin{equation}\label{9}
v^2=\frac{(-2h\sum_{i=1}^n a_iD_i^5)}{(\sum_{i=1}^m b_iB_iC_i^2)} u^4-\frac{(20h^3\sum_{i=1}^n a_iD_i^3+ \sum_{i=1}^m b_iB_i^3)}{(3\sum_{i=1}^m b_iB_iC_i^2)}u^2-\frac{(2h^5\sum_{i=1}^n a_iD_i)
}{(\sum_{i=1}^m b_iB_iC_i^2)}.
\end{equation}
\\

We see that if
\begin{equation*}
L':=\frac{(-2h^5\sum_{i=1}^n a_iD_i)
}{(\sum_{i=1}^m b_iB_iC_i^2)}
\end{equation*}
\\

 to be square, say $q^2$, (This is possible by fixing the other parameters and  choosing appropriate values for $h$.),
we may transform (similar to the previous theorem) the above elliptic curve to an elliptic curve of  the form

 $y^2+a_1xy+a_3y=x^3+a_2x^2+a_4x+a_6$, where $a_i\in \mathbb{Q}$.
\\

If the above elliptic curve has positive rank (For every arbitrary nonzero integers $a_i$ and $b_i$, this is done  by choosing appropriate values for $B_i$, $C_i$, $D_i$, $h$.), by calculating $x_i$, $y_i$, $X_i$, $Y_i$, from the relations $x_i=D_iu+h$, $y_i=B_iu-C_iv$, $X_i=-D_iu+h$, and $Y_i=B_iu+C_iv$, and after  some simplifications  and canceling  the  denominators of $x_i=D_iu+h$, $y_i=B_iu-C_iv$, $X_i=-D_iu+h$, and $Y_i=B_iu+C_iv$, we may obtain infinitely many integer solutions for the  Diophantine equation. Now the proof is complete.
\end{theorem}

Also we know that if
 $(x_1, \cdots,x_n,y_1, \cdots ,y_m,X_1, \cdots, X_n,Y_1, \cdots, Y_m)$
is a solution for the Diophantine equations, then for every arbitrary $t$,

 $(tx_1, \cdots,tx_n,t^2y_1, \cdots ,t^2y_m,tX_1, \cdots, tX_n,t^2Y_1, \cdots,t^2 Y_m)$
 is a solution too.

Then if we obtain a rational solution, we may get an integer solution for  the Diophantine equations by multiplying the both sides of the Diophantine equations by an appropriate $t$.
\\

Finally we mention that each point on the elliptic curve can be represented in the form $(\frac{r}{s^2},\frac{t}{s^3})$, where $r$, $s$, $t$ $\in \mathbb{Z}$.\\
So if we put
$nP=(\frac{r_n}{s_n^2},\frac{t_n}{s_n^3})$, that the point P is one of  the elliptic curve generators, we may  obtain a parametric solution for each case of Diophantine equations by using this new point $P'=nP=(\frac{r_n}{s_n^2},\frac{t_n}{s_n^3})$. Also by choosing  the other appropriate values of $B_i$, $C_i$, $D_i$, $h$, and getting the new elliptic curve of rank $\geq 1$( and repeating   the above process) , we may obtain infinitely many nontrivial parametric solutions for each case of  the above Diophantine equations.

We use the Sage software for calculating  the rank of the elliptic curves. (see \cite{4})


\begin{thebibliography}{HD}
\bibitem[1] {1}  L.  E.  DICSON,  History  of the  Theory of Numbers, Vol. II:  Diophantine  Analysis,  G.  E.  STECHERT.
  Co., New York,  $1934$.
\bibitem[2]{2} N. ELKIES, $ "On  A4 + B4 + C4 = D4"$. Mathematics of Computation. $1988$. $51 (184): 825–835.$
\bibitem[3]{3}  L.  J.  LANDER  and  T.  R.  PARKIN,  "Counterexamples  to  Euler's  conjecture  on  sums  of  like
powers,"  Bull.  Amer.  Math.  Soc, $1966,  VOL. 72,   p.  1079$.
\bibitem[4]{4} SAGE software, available from http://sagemath.org.
\bibitem[5] {5}  L. C.  WASHINGTON, Elliptic Curves: Number Theory and Cryptography, Chapman-Hall, $2008$.







\normalsize
\baselineskip=17pt



\bibitem[]{}

\bibitem[]{}

\bibitem[]{}

\end{thebibliography}
\end{document}